\documentclass[a4paper,12pt]{article}
\title{Weierstrass semigroups on double covers of genus four curves}
\author{Seon Jeong Kim \thanks{E-mail: skim@gnu.ac.kr}\\
{\small Department of Mathematics and RINS, Gyeongsang National University}\\
{\small  Jinju, 660-701, Republic of Korea}\\
and\\
Jiryo Komeda\thanks{E-mail: komeda@gen.kanagawa-it.ac.jp}\\
{\small Department of Mathematics, Kanagawa Institute of Technology}\\
{\small  Atsugi, 243-0292, Japan}}
\date{}

\usepackage{amssymb}
\usepackage[leqno]{amsmath}
\usepackage{latexsym}
\usepackage{theorem}
\usepackage{enumerate}
\newtheorem{theorem}{Theorem}[section]
\newtheorem{corollary}[theorem]{Corollary}
\newtheorem{proposition}[theorem]{Proposition}
\newtheorem{lemma}[theorem]{Lemma}

\newtheorem{remark}[theorem]{Remark}

\newcounter{examplec}[section]

\newcommand{\qed}{\hfill $  \Box $}

\newcommand{\tH}{\tilde{H}}
\newcommand{\tC}{\tilde{C}}
\newcommand{\tP}{\tilde{P}}
\newcommand{\tg}{\tilde{g}}
\newcommand{\tih}{\tilde{h}}
\newcommand{\la}{\langle}
\newcommand{\ra}{\rangle}
\newcommand{\pr}{\vskip3mm
\noindent
{\it Proof.} }

\newcommand{\smax}{s_{{\rm max}}}

\newcommand{\dis}{\displaystyle}
\newcommand{\fn}{\frac{n-1}{2}}
\newcommand{\NI}{\mathbb{N}_0}

\newcommand{\DC}{\mbox{{\it of double covering type}}}
\newcommand{\dc}{\mbox{of double covering type}}
\newcommand{\com}{\mbox{, }}
\Roman{equation}
\begin{document}
\maketitle
\renewcommand{\thefootnote}{\fnsymbol{footnote}}
\footnotetext{The first author is partially supported by Basic Science Research Program through the
National Research Foundation of Korea(NRF) funded by the Ministry of
Education, Science and Technology (2011-0027246).
The second author is partially supported by
Grant-in-Aid for Scientific Research (24540057), Japan Society for the Promotion Science. }

\renewcommand{\abstractname}{}
\begin{abstract}

Let $C$ be a complete non-singular irreducible curve of genus $4$ over an algebraically closed field of characteristic $0$.
We determine all possible Weierstrass semigroups of ramification points on double covers of $C$ which have genus greater than $11$.
Moreover, we construct double covers with ramification points whose Weierstrass semigroups are the possible ones.

\vspace{2mm} \noindent {\bf 2010 Mathematics Subject
Classification:} 14H55, 14H30, 14H45 \\
{\bf Key words:} Weierstrass semigroup of a point, Double cover of a curve, Curve of genus $4$ \end{abstract}
\section{Introduction} 
\label{intro}

A submonoid $H$ of the additive monoid $\mathbb{N}_0$ of non-negative integers is called a {\it numerical semigroup} if its complement $\mathbb{N}_0\backslash H$ is a finite set.
The cardinality of $\mathbb{N}_0\backslash H$ is called the {\it genus} of $H$, which is denoted by $g(H)$.
Let $C$ be a complete nonsingular irreducible curve over an algebraically closed field $k$ of characteristic 0, which is called a {\it
curve} in this paper.
Let $k(C)$ be the field of rational functions on
$C$.
For a point $P$ of $C$, we set
$$H(P)=\{\alpha\in \mathbb{N}_{0}|\mbox{ there exists } f\in k(C)\mbox{ with } (f)_{\infty}=\alpha P\},$$
which is called the {\it Weierstrass semigroup of $P$}.
It is known that the Weierstrass semigroup of a point on a curve of genus $g$ is a numerical semigroup of genus $g$.

Let $\pi:{\tilde C}\longrightarrow C$ be a double covering of a curve.
We are interested in the Weierstrass semigroup $H({\tilde P})$ of a ramification point ${\tilde P}$ on the double cover $\tC$ of $C$.
A numerical semigroup is said to be {\it of double covering type}, if it is such a Weierstrass semigroup $H({\tilde P})$. For a numerical semigroup $\tH$ we denote by $d_2(\tH)$ the set consisting of the elements  $\tih/2$ with even $\tih\in \tH$, which becomes a numerical semigroup.
Using this notation we have $d_2(H({\tilde P}))=H(\pi({\tilde P}))$.
Let $\tg$ be the genus of ${\tilde C}$.
If $C$ is the projective line, then ${\tilde P}$ is a Weierstrass point on the hyperelliptic curve ${\tilde C}$ for any $\tg\geqq 2$.
Hence, the semigroup $H({\tilde P})$ is generated by $2$ and $2\tg+1$.
If $C$ is an elliptic curve, then the semigroup $H({\tilde P})$ is either $\la 3,4,5\ra$ or $\la 3,4\ra$ or $\la 4,5,6,7\ra$ or $\langle 4,6,2\tg-3\rangle$ with $\tg\geqq 4$ or $\langle 4,6,2\tg-1,2\tg+1\rangle$ with $\tg\geqq 4$.
Here for any positive integers $a_{1},a_{2},\ldots,a_{n}$ we denote by $\langle a_{1},a_{2},\ldots,a_{n}\rangle$ the additive monoid generated by $a_{1},a_{2},\ldots,a_{n}$.
Conversely, there is a double cover of an elliptic curve with a ramification point whose Weierstrass semigroup is any semigroup of the above ones (for example, see ~\cite{kom1}, ~\cite{kom2}).
Oliveira and Pimentel ~\cite{oli-pim} studied Weierstrass semigroups $H({\tilde P})$ in the case where the genus of $C$ is $2$.
They showed that for a semigroup $\tH=\langle 6,8,10,n\rangle$ with an odd number $n\geqq 11$ there exists a double covering $\pi:{\tilde C}\longrightarrow C$ with a ramification point ${\tilde P}$ such that $H({\tilde P})=\tH$.
Moreover, in ~\cite{kom2} we showed that all numerical semigroups $\tH$ of genus $\geqq 6$ satisfying $g(d_2(\tH))=2$ are \dc.
In ~\cite{kom4} we proved that every numerical semigroup $\tH$ of genus $\geqq 9$  with $g(d_2(\tH))=3$ is \dc.

In this paper we will study the Weierstrass semigroups of ramification points on double covers of genus four curves.
Namely we prove 
\vskip3mm
{\bf Main Theorem.} {\it Let $H$ be a numerical semigroup of genus $4$ and $\tH$ a numerical semigroup with $d_2(\tH)=H$.
If $g(\tH)\geqq 12$, then $\tH$ is $\DC$.}
\vskip3mm
But we note that Main Theorem does not hold in the case where $g\geqq 5$.
That is to say, for any $g\geqq 5$ there are numerical semigroups $\tH$ with $g(d_2(\tH))=g$ and $g(\tH)\geqq 3g$ which are not \dc \ (see ~\cite{kom3}).
\vskip3mm

In Section 2 we will treat known facts and new results which work well for numerical semigroups of any genus with some properties. 
We know that a numerical semigroup of genus $4$ is either $\la 2,9\ra$ or
$\la 3,5 \ra$ or $\la 3,7,8 \ra$ or $\la 4,5,6 \ra$ or $\la 4,5,7 \ra$ or $\la 4,6,7,9 \ra$ or $\la 5,6,7,8,9 \ra$.
By ~\cite{kom-ohb2} any numerical semigroup $\tH$ with $d_2(\tH)=\la 2,9\ra$ is \dc.
From Section 3 to Section 8 we will prove Main Theorem for each numerical semigroup $\tH$ with $g(d_2(\tH))=4$ and $d_2(\tH)\not=\la 2,9\ra$.
\section{General theory} 
In this section we will state results in \cite{kom-ohb1} which are used after Section 2 and give criterions for a numerical semigroup with some conditions to be of double covering type.
To explain these we need some terminology and notation.
For a numerical semigroup $H$ we denote by $c(H)$ and $n(H)$ the least integer $c$ with $c+\NI\subseteqq H$ and the least odd integer in $H$ respectively.
Let $m$ be a positive integer.
An {\it $m$-semigroup} means a numerical semigroup whose least positive integer is $m$.
For an $m$-semigroup $H$ we set
$S(H)=\{m,s_1,\ldots,s_{m-1}\}$
where $s_i=\min\{h\in H\mid h\equiv i\mbox{ mod }m\}$ for $i=1,\ldots,m-1$.
We call $S(H)$ the {\it standard basis} for $H$.
Moreover, we set $\smax=\max\{s_1,\ldots,s_{m-1}\}$.
In Lemma 2.1 of ~\cite{kom-ohb1} we may replace the assumption $n> 2c(H)-2$ by $n>c(H)+m-2$ as follows:
\begin{remark} \label{SBasis} 
Let $H$ be an $m$-semigroup and $n$ an odd integer larger than $c(H)+m-2$.
We set
$$H_{n}=2H+n\mathbb{N}_{0}=\{2h+ni\mid h\in H\com i\in \NI\}.$$
Assume that $n> 2m$, which holds for any $H\not=\la m,m+1,\ldots,m+m-1\ra$.
Then we have the following:\\
{\rm \romannumeral1)} $H_{n}$ is a $2m$-semigroup with the standard basis
$$S(H_{n})=\{2m,2s_{1},\ldots,2s_{m-1},n,n+2s_{1},\ldots,n+2s_{m-1}\}.$$
{\rm \romannumeral2)} The genus of $H_{n}$ is $\displaystyle 2g(H)+\frac{n-1}{2}$.
\end{remark}
\vskip3mm
Using Remark 2.1 we get the following whose proof is given in Lemma 2.3 of ~\cite{kom4}:
\begin{lemma}\label{ineq} 
Let $H$ be an $m$-semigroup and $\tH$ a numerical semigroup with $d_2(\tH)=H$.
If $n=n(\tH)\geqq c(H)+m-1$ and $n\not=2m-1$, then we get
$$g(H)+\frac{n-1}{2}\leqq g(\tH)\leqq 2g(H)+\frac{n-1}{2}.$$
The right-hand $($resp. left-hand$)$ equality of the above inequalities holds if and only if $\tH=2H+n\NI$  $($resp. $2H+n\NI+(n+2)\NI+(n+4)\NI+\cdots+(n+2(m-1))\NI\mbox{})$.
\end{lemma}
In Theorem 2.2 of ~\cite{kom-ohb1} we may replace the assumption $n\geqq 2c(H)-1$ by $n\geqq \max\{c(H)+m-1,2r+1\}$ where $r$ is the genus of $H$ as follows:
\begin{proposition} \label{(n-1)/2} 
Let $H$ be a Weierstrass $m$-semigroup of genus $r\geqq 0$.
Let $n$ be an odd number with $n\geqq \max\{c(H)+m-1,2r+1\}$.
Then the numerical semigroup $H_{n}=2H+n\mathbb{N}_{0}$ is \DC.
\end{proposition}

We get the following statement for special numerical semigroups $\tH$ with $\displaystyle g(\tH)=2g(H)+\frac{n-1}{2}-1$ in ~\cite{kom3}.
\begin{proposition} \label{smax-m} 
 Let $H$ be an $m$-semigroup and $n$ an odd integer with $n\geqq \max\{c(H)+m-1,2g(H)+3\}$.
 We set $\tH=2H+n\NI+(n+2(\smax-m))\NI$.
 Then we have the following:\\
 {\rm 1)} $\displaystyle g(\tH)=2g(H)+\frac{n-1}{2}-1$.\\
 {\rm 2)} If $H$ is Weierstrass, then $\tH$ is \dc.
 \end{proposition}

If the genus of $\tH$ is $\displaystyle 2g(H)+\frac{n-1}{2}-1$ and $d_2(\tH)$  is symmetric, then we can solve our problem (see ~\cite{kom3}):
\begin{proposition} \label{symm-1} 
Let $H$ be a symmetric $m$-semigroup and $\tH$ a numerical semigroup with $d_2(\tH)=H$.
Assume that $n=n(\tH)\geqq \max\{2g(H)+m-1,2g(H)+3\}$.
Then the following are equivalent: \\
{\rm 1)} $\displaystyle g(\tH)=2g(H)+\frac{n-1}{2}-1$,
\\
{\rm 2)} $\tH=2H+n\mathbb{N}_0+(n+2(\smax-m))\mathbb{N}_0.$ \\
In this case, if $H$ is Weierstrass, then $\tH$ is \dc.
\end{proposition}

We consider the special case where the genus of $\tH$ is $\displaystyle 2g(H)+\frac{n-1}{2}-2$ and $d_2(\tH)$ is symmetric.
Then we get the following:
\begin{proposition} \label{symm-2} 
Let $H$ be a symmetric $m$-semigroup with $m\geqq 3$.
Let $n$ be an odd integer satisfying $n\geqq 2g(H)+m-1$.
We set
$$\tH=2H+n\NI+(n+2(\smax-2m))\NI.$$
Then we have the following: \\
{\rm 1)} The genus of $\tH$ is $\dis 2g(H)+\frac{n-1}{2}-2$.\\
{\rm 2)} If $H$ is Weierstrass and $g(\tH)\geqq 3g(H)$, then $\tH$ is \DC.
\end{proposition}
\pr
1) Let $s_i\not=\smax$ and $n+2s_i-2m=n+2(\smax-2m)+\tih$ with $\tih\in \tH$.
Then we get $2s_i=2\smax-2m+\tih$.
By Proposition 3.4 and Lemma 3.2 in ~\cite{kom3} there exists $s_j\not=\smax$ such that $2\smax=2s_i+2s_j=2\smax-2m+\tih+2s_j$.
Hence, we get $2m=\tih+2s_j\geqq 2s_j>2m$, which is a contradiction.
Therefore, we have $\dis g(\tH)=2g(H)+\frac{n-1}{2}-2$.

2) Let $(C,P)$ be a pointed curve of genus $g$ with $H(P)=H$.
We set $\dis D=\frac{n+1}{2}P-Q_1-Q_2$, where $Q_i\not=P$ for $i=1,2$.
Then we have $h^0(mP-Q_1-Q_2)=0$ for almost all points $Q_2$.
Since $\dis 2g+\frac{n-1}{2}-2=g(\tH)\geqq 3g$, we obtain $n\geqq 2g+5$.
Hence, we have $\deg(2D-P)=n-4\geqq 2g+1$.
Thus, the divisor $2D-P$ is very ample.
Hence, we can construct a double covering $\pi:\tC\longrightarrow C$ with a ramification point $\tP$ over $P$.
Since the degree of the branch divisor of $\pi$ is $\deg (2D)$, the genus of $\tC$ is equal to $g(\tH)$.
We have
$$h^0((n-1)\tP)=h^0(\frac{n-1}{2}P)+h^0(Q_1+Q_2-P)$$
and
$$h^0((n+1)\tP)=h^0(\frac{n+1}{2}P)+h^0(Q_1+Q_2).$$
In view of $Q_i\not=P$ we have $n\in H(\tP)$.
Moreover, we obtain
$$h^0((n-1+2(\smax-2m))\tP)=h^0((\frac{n-1}{2}+\smax -2m)P)+h^0((\smax-2m-1)P+Q_1+Q_2)$$
and
$$h^0((n+1+2(\smax-2m))\tP)=h^0((\frac{n+1}{2}+\smax -2m)P)+h^0((\smax-2m)P+Q_1+Q_2).$$
On the other hand, we have
$$h^0((\smax-2m-1)P+Q_1+Q_2)$$
$$=\smax-g-2m+2+h^0(mP-Q_1-Q_2)=\smax-g-2m+2$$
and
$$h^0((\smax-2m)P+Q_1+Q_2)$$
$$=\smax-g-2m+3+h^0((m-1)P-Q_1-Q_2)=\smax-g-2m+3$$
because of $K\sim (\smax-m-1)P$ and $h^0(mP-Q_1-Q_2)=0$.
Hence, we get $n+2(\smax-2m)\in H(\tP)$.
Since the genus of $H(\tP)$ coincides with that of $\tH$, we have $\tH=H(\tP)$, which implies that $\tH$ is \dc.
\qed
\vskip3mm

The following is a key tool for showing that a given numerical semigroup $\tH$ with $g(d_2(\tH))=4$ and $g(\tH)\geqq 3g(d_2(\tH))$ is \dc.
\begin{theorem}\label{w3gdc} 
Let $H$ be a Weierstrass numerical semigroup.
Take a pointed curve $(C,P)$ with $H(P)=H$.
Let $\tH$ be a numerical semigroup with $d_2(\tH)=H$ with $g(\tH)\geqq 3g(H)$.
Let $r$ be an integer with $0\leqq r\leqq g(H)$ and $Q_1,\ldots,Q_r$ points of $C$ different from $P$ with $h^0(Q_1+\cdots+Q_r)=1$.
We set $n=n(\tH)$.
Assume that $\dis g(\tH)=2g(H)+\frac{n-1}{2}-r$ and that $\tH$ has an expression
$$\tH=2H+\la n,n+2l_1,\ldots,n+2l_s\ra$$
of generators with positive integers $l_1,\ldots,l_s$ such that
$$h^0(l_iP+Q_1+\cdots+Q_r)=h^0((l_i-1)P+Q_1+\cdots+Q_r)+1$$
for all $i$.
Then there is a double covering $\pi:\tC\longrightarrow C$ with a ramification point $\tP$ over $P$ satisfying $H(\tP)=\tH$, hence $\tH$ is $\DC$.
\end{theorem}
\pr
By Theorem 2.2 in ~\cite{kom4} it suffices to show that the complete linear system $|nP-2Q_1-\cdots-2Q_r|$ is base point free.
Since $g(\tH)\geqq 3g(H)$ and $\displaystyle g(\tH)=2g(H)+\frac{n-1}{2}-r$, we have $n-2r\geqq 2g(H)+1$.
Hence, the complete linear system is very ample, which implies that it is base point free.
\qed
\vskip3mm
Using the above theorem we get the following:
\begin{corollary}\label{lowest}
Let $H$ be a Weierstrass numerical semigroup and $\tH$ a numerical semigroup with $d_2(\tH)=H$ and $g(\tH)\geqq 3g(H)$.
Assume that $n=n(\tH)\geqq c(H)+m-1$ with $n\not=2m-1$.
If $\dis g(\tH)=g(H)+\frac{n-1}{2}$, then $\tH$ is \dc.
\end{corollary}
\pr
Let $(C,P)$ be a pointed curve with $H(P)=H$.
Let $Q_1,\ldots,Q_g$ be general points of $C$.
Then we get $h^0(Q_1+\cdots+Q_g)=1$.
In view of $\dis g(\tH)=g(H)+\frac{n-1}{2}$, we have
$$\tH= 2H+n\NI+(n+2)\NI+(n+4)\NI+\cdots+(n+2(m-1))\NI$$
by Lemma \ref{ineq}.
Since we obtain $h^0(iP+Q_1+\cdots+Q_g)=i+1$ for $i\geqq 0$, by Theorem \ref{w3gdc} $\tH$ is \dc.
\qed
\section{Numerical semigroups $\tH$ with $d_2(\tH)=\la 3,5\ra$
} 
\begin{theorem} 
We set $H=\la 3,5\ra$.
Let $\tH$ be a numerical semigroup with $d_2(\tH)=H$.
If $g(\tH)\geqq 12$, then $\tH$ is \dc.
\end{theorem}
\pr
We take a pointed curve $(C,P)$ with $H(P)=\la 3,5\ra$.
Then a canonical divisor $K$ on $C$ is linearly equivalent to $6P$.
We set $n=n(\tH)$.
If $n\leqq 7$, then we have $g(\tH)<12$.
If $n=9$, then $\tH\supseteqq \la 6,9,10\ra$.
We have $g(\la 6,9,10\ra)=12$.
If we set $D=5P$, the same proof as in Theorem 2.2 of ~\cite{kom-ohb1} works well.
Hence, the semigroup $\tH$ is \dc.
Therefore, we may assume  that $n\geqq 11$.
By Lemma \ref{ineq} we obtain $\dis 4+\frac{n-1}{2}\leqq g(\tH)\leqq 8+\frac{n-1}{2}$.

If $\dis g(\tH)=8+\frac{n-1}{2}$,
then $\tH$ is \dc \ by Proposition \ref{(n-1)/2}.

If $\dis g(\tH)=7+\frac{n-1}{2}$, then $\tH$ is \dc \ by Proposition \ref{symm-1}, because the semigroup $H$ is symmetric.

Let $\dis g(\tH)=6+\frac{n-1}{2}$.
First, we note that
$$S(2H+n\NI)=\{6,10,20,n,n+10,n+20\}.$$
Hence, $\tH$ is either
$$2H+n\NI+(n+8)\NI \mbox{ or } 2H+n\NI+(n+4)\NI.$$
In the first case, $\tH$ is \dc \ by Proposition \ref{symm-2},
because $H$ is symmetric and $8$ is $2(s_{max}-2\cdot 3)$.
We consider the second case. 
Let $S$ be a point distinct from $P$.
Then there exists a unique set $\{Q_1,Q_2\}$ of points of $C$ distinct from $P$ such that $3P\sim S+Q_1+Q_2$.
Then we have
$$h^0(P+Q_1+Q_2)=h^0(K-P-Q_1-Q_2)=h^0(2P+3P-Q_1-Q_2)=h^0(2P+S)=1,$$
because $C$ has a unique trigonal system $|3P|$.
Moreover, we have
$$h^0(2P+Q_1+Q_2)=1+h^0(P+3P-Q_1-Q_2)=1+h^0(P+S)=2.$$
Hence, by Theorem \ref{w3gdc} the semigroup $2H+n\NI+(n+4)\NI$ is of double covering type.

Let $\dis g(\tH)=5+\frac{n-1}{2}$.
Then $\tH$ is either $\la 6,10,n,n+4,n+8\ra$ or $\la 6,10,n,n+2\ra$.
In the first case we will show that $\tH$ is \dc.
Let $Q_1$ be a point of $C$ distinct from $P$.
If $Q_2$ and $Q_3$ are general points, then we have
$$h^0(5P-Q_1-Q_2-Q_3)=0.$$
Hence, we get 
$$h^0(2P+Q_1+Q_2+Q_3)=2+h^0(4P-Q_1-Q_2-Q_3)=2.$$
Moreover, we have
$$h^0(P+Q_1+Q_2+Q_3)=1+h^0(5P-Q_1-Q_2-Q_3)=1,$$
which implies that
$$h^0(2P+Q_1+Q_2+Q_3)=h^0(P+Q_1+Q_2+Q_3)+1.$$
We also see that
$$h^0(3P+Q_1+Q_2+Q_3)=3+h^0(3P-Q_1-Q_2-Q_3)=3.$$
Moreover, we get $h^0(4P+Q_1+Q_2+Q_3)=4$, which implies that
$$h^0(4P+Q_1+Q_2+Q_3)=h^0(3P+Q_1+Q_2+Q_3)+1.$$
By Theorem \ref{w3gdc} the semigroup $\la 6,10,n,n+4,n+8\ra$ is \dc.
Next, we will consider the second case.
If $Q_1$ is distinct from $P$ and $Q_2$ is a general point, then we have $h^0(5P-Q_1-Q_2)=1$, which implies that an effective divisor $F$ with $F\sim 5P-Q_1-Q_2$ is uniquely determined.
We set $F=R_1+R_2+R_3$.
Then we have $R_1+R_2+R_3\not=3P$.
Hence, we may assume that $R_1\not=P$.
We set $Q_3=R_1$.
Then we have
$$h^0(P+Q_1+Q_2+Q_3)=1+h^0(5P-Q_1-Q_2-Q_3)=1+h^0(R_2+R_3)=2.$$
Moreover, we have $h^0(Q_1+Q_2+Q_3)=1$ for general points $Q_1$ and $Q_2$.
In fact, we assume that $h^0(Q_1+Q_2+Q_3)=2$ for general points $Q_1$ and $Q_2$.
Then we have $Q_1+Q_2+Q_3\sim 3P$, which contradicts $h^0(3P)=2$.
Hence, we get
$$h^0(P+Q_1+Q_2+Q_3)=h^0(Q_1+Q_2+Q_3)+1.$$
By Theorem \ref{w3gdc} the semigroup $\la 6,10,n,n+2\ra$ is \dc.

Let $\dis g(\tH)=4+\frac{n-1}{2}$.
By Corollary \ref{lowest} $\tH$ is \dc.
\qed
\section{Numerical semigroups $\tH$ with $d_2(\tH)=\la 3,7,8\ra$} 
\begin{theorem} 
We set $H=\la 3,7,8\ra$.
Let $\tH$ be a numerical semigroup with $d_2(\tH)=H$.
If $g(\tH)\geqq 12$, then $\tH$ is \dc.
\end{theorem}
\pr
We set $n=n(\tH)$.
Let $n\leqq 7$.
Then $\tH$ contains either $\la 6,7,16\ra$ or $\la 5,6,14\ra$ or $\la 3,14,16\ra$, whose genus is less than $12$.
Hence, we may assume $n\geqq 9$.
Let $(C,P)$ be a pointed curve with $H(P)=\la 3,7,8\ra$.
Then we have $K_C\sim 4P+R_1+R_2$ where $R_1$ and $R_2$ are points distinct from $P$.

Let $\dis g(\tH)=8+\frac{n-1}{2}$.
By Lemma \ref{ineq} and Proposition \ref{(n-1)/2}  $\tH$ is of double covering type.

Let $\dis g(\tH)=7+\frac{n-1}{2}$.
First, we note that
$$S(2H+n\NI)=\{6,14,16,n,n+14,n+16\}.$$
Then $\tH$ is either
$2H+n\NI+(n+10)\NI \mbox{ or } 2H+n\NI+(n+8)\NI.$
By Proposition \ref{smax-m} the first one is \dc.
Let $Q=R_1$.
Then we have
$$h^0(3P+Q)=1+h^0(P+R_2)=2\mbox{ and }h^0(4P+Q)=2+h^0(R_2)=3,$$
which implies that $h^0(4P+Q)=h^0(3P+Q)+1$.
Hence, by Theorem \ref{w3gdc} the latter one is also \dc.

Let $\dis g(\tH)=6+\frac{n-1}{2}$.
$\tH$ is either $2H+n\NI+(n+2)\NI$ or $2H+n\NI+(n+4)\NI$ or $2H+n\NI+(n+8)\NI+(n+10)\NI$.
Let $Q_1=R_1$ and $Q_2=R_2$.
Then we have
$$h^0(P+Q_1+Q_2)=h^0(3P)=2=h^0(Q_1+Q_2)+1,$$
which implies that the first semigroup is \dc.
Let $Q_1$ be a distinct point from $P$, $R_1$ and $R_2$.
Then we will show that $h^0(2P+R_1+R_2-Q_1)=1$.
In fact, if $h^0(2P+R_1+R_2-Q_1)=2$, we would have either
$$2P+R_1+R_2-Q_1\sim 3P\mbox{ or }2P+R_1+R_2-Q_1\sim P+R_1+R_2.$$
The first case implies that $R_1+R_2\sim P+Q_1$, which is a contradiction.
The second case means that $P\sim Q_1$.
This contradicts $P\not=Q_1$.
Hence, there is a unique effective divisor $E$ with $E\sim 2P+R_1+R_2-Q_1$.
Then $E$ contains a point $S$ different from $P$.
If $E=3P$, then $3P\sim 2P+R_1+R_2-Q_1$, which implies that $P+Q_1\sim R_1+R_2$.
This is a contradiction.
For any point $Q_2$ we have
$h^0(3P+R_1+R_2-Q_1-Q_2)=1$ .
As a matter of fact, if there is some point $Q_2$ such that $h^0(3P+R_1+R_2-Q_1-Q_2)=2$, then either $3P+R_1+R_2-Q_1-Q_2\sim 3P$ or $3P+R_1+R_2-Q_1-Q_2\sim P+R_1+R_2$ holds.
This contradicts the assumption on $Q_1$.
Let $Q_2=S$.
Then we obtain
$$h^0(P+Q_1+Q_2)=3+1-4+h^0(3P+R_1+R_2-Q_1-Q_2)=1$$
and
$$h^0(2P+Q_1+Q_2)=1+h^0(2P+R_1+R_2-Q_1-Q_2)=1+h^0(E-S)=2.$$
Hence, the second semigroup is \dc.
We consider the third case.
Let $Q_1$ be distinct from $P$, $R_1$ and $R_2$.
Since we have $h^0(P+R_1+R_2)=2$, there is a unique effective divisor $E$ which is linearly equivalent to $P+R_1+R_2-Q_1$.
Let $Q_2$ be a point different from $P$ with $Q_2\not< E$.
Then we have
$h^0(5P+Q_1+Q_2)=7+1-4=4$
and
$$h^0(4P+Q_1+Q_2)=6+1-4+h^0(R_1+R_2-Q_1-Q_2)=3,$$
because $Q_1$ is distinct from $R_1$ and $R_2$.
Moreover, we obtain
$$h^0(3P+Q_1+Q_2)=2+h^0(P+R_1+R_2-Q_1-Q_2)=2+h^0(E-Q_2)=2,$$
because of $Q_2\not< E$.
Hence, the third semigroup is \dc.

Let $\dis g(\tH)=5+\frac{n-1}{2}$.
Then the semigroup $\tH$ is either
$$2H+n\NI+(n+2)\NI+(n+10)\NI\mbox{ or }
2H+n\NI+(n+4)\NI+(n+8)\NI.$$
Let $Q_1$, $Q_2$ and $Q_3$ be points distinct from $P$ where $Q_1$ and $Q_2$ are general.
Then we have
$$h^0(5P+Q_1+Q_2+Q_3)=5\mbox{ and }h^0(4P+Q_1+Q_2+Q_3)=4.$$
Since $Q_1$ and $Q_2$ are general, we get
$h^0(K-P-Q_1-Q_2)=1$, which implies that
there is a unique effective divisor $E$ with $E\sim K-P-Q_1-Q_2$ and $E\not=3P$.
Let $Q_3<E$.
Then we obtain
$$h^0(P+Q_1+Q_2+Q_3)=1+h^0(K-P-Q_1-Q_2-Q_3)=1+h^0(E-Q_3)=2.$$
Moreover, we have
$h^0(Q_1+Q_2+Q_3)=1$,
because $Q_1$ and $Q_2$ are general.
Hence, the first one is \dc.
Let $Q_1$, $Q_2$ and $Q_3$ be general points.
Then we have
$$h^0(4P+Q_1+Q_2+Q_3)=4\mbox{ and }h^0(3P+Q_1+Q_2+Q_3)=3.$$
Moreover, we obtain
$$h^0(2P+Q_1+Q_2+Q_3)=2\mbox{ and }h^0(P+Q_1+Q_2+Q_3)=1.$$
Thus, the latter one is also of double covering type.

Let $\dis g(\tH)=4+\frac{n-1}{2}$.
By Corollary \ref{lowest} $\tH$ is \dc.
\qed
\section{Numerical semigroups $\tH$ with $d_2(\tH)=\la 4,5,6\ra$} 
\begin{theorem} 
We set $H=\la 4,5,6\ra$.
Let $\tH$ be a numerical semigroup with $d_2(\tH)=H$.
If $g(\tH)\geqq 12$, then $\tH$ is \dc.
\end{theorem}
\pr
We take a pointed curve $(C,P)$ with $H(P)=\la 4,5,6\ra$.
In this case we have $K\sim 6P$.
If $n=n(\tH)\leqq 9$, then $\tH\supseteqq \la 8,9,10,12\ra$.
Since the genus of  $\la 8,9,10,12\ra$ is $12$, we have $\tH=\la 8,9,10,12\ra$.
If we set $D=5P$, the same proof as in Theorem 2.2 of ~\cite{kom-ohb1} works well.
Hence, the semigroup $\tH$ is \dc.
If $n\leqq 7$, then we have $g(\tH)<12$.
Thus, we may assume  that $n\geqq 11$.

If $\dis g(\tH)=8+\frac{n-1}{2}$ or $\dis 7+\frac{n-1}{2}$, then $\tH$ is of double covering type by Propositions \ref{(n-1)/2} and \ref{symm-1}.

If $\dis g(\tH)=6+\frac{n-1}{2}$, then $\tH$ is either
$$2H+n\NI+(n+2)\NI \mbox{ or } 2H+n\NI+(n+4)\NI \mbox{ or } 2H+n\NI+(n+6)\NI,$$
because of
$S(2H+n\NI)=\{8,10,12,22,n,n+10,n+12,n+22\}$.
By Proposition \ref{symm-2} the last semigroup is \dc.
We will show that the first semigroup is \dc.
There exists a point $Q_1\in C$ such that $Bs|5P-Q_1|\not=\emptyset.$
In fact, assume that we had $Bs|5P-Q_1|=\emptyset$ for any point $Q_1\in C$.
Then $|5P|$ is a very ample divisor with $\dim |5P|=2$.
But $C$ is not isomorphic to a plane curve.
This is a contradiction.
Take a point $Q_1\in C$ such that $Bs|5P-Q_1|\not=\emptyset.$
Then $Q_1\not=P$, because $|4P|$ is base point free.
Hence, there is a point $Q_2\in C$ such that $\dim |5P-Q_1-Q_2|=\dim |5P-Q_1|$.
In this case, we have $Q_2\not=P$, because $|4P|$ and $|5P|$ are base point free.
Then we obtain
$$h^0(P+Q_1+Q_2)=h^0(5P-Q_1-Q_2)=h^0(5P-Q_1)=2$$
and $h^0(Q_1+Q_2)=1$.
Hence, the first one is \dc.
We will show that the second semigroup is \dc.
For any point $Q\not=P$, the set consisting of three points $R_1(Q)$, $R_2(Q)$ and $R_3(Q)$ with $4P-Q\sim R_1(Q)+R_2(Q)+R_3(Q)$ is uniquely determined, where $R_1(Q)$, $R_2(Q)$ and $R_3(Q)$ are distinct from $P$ because of $h^0(3P)=1$.
Then we have
$$h^0(4P-Q-R_1(Q))=h^0(R_2(Q)+R_3(Q))=1.$$
Moreover, we have
$$h^0(5P-Q-R_1(Q))=h^0(P+R_2(Q)+R_3(Q)).$$
Assume that for any point $Q\not=P$ we had $h^0(5P-Q-R_1(Q))=2$.
Since we have at most two trigonal systems on $C$, there exists a  point $Q_1\not=P$ such that for infinitely many points $Q'$ we have 
$$P+R_2(Q_1)+R_3(Q_1)\sim P+R_2(Q')+R_3(Q'),$$
which implies that $R_2(Q_1)+R_3(Q_1)\sim R_2(Q')+R_3(Q')$.
Therefore, we get
$Q_1+R_1(Q_1)\sim Q'+R_1(Q')$
for infinitely many points $Q'$.
Hence, we get $Q_1=R_1(Q')$ and $R_1(Q_1)=Q'$.
But the set $\{R_1(Q_1),R_2(Q_1),R_3(Q_1)\}$ is uniquely determined.
This is a contradiction.
Thus, there exists a point $Q\not=P$ with $h^0(5P-Q-R_1(Q))=1$ where $R_1(Q)\not=P$.
Hence, we get
$$h^0(2P+Q+R_1(Q))=1+h^0(4P-Q-R_1(Q))=2$$
and
$$h^0(P+Q+R_1(Q))=h^0(5P-Q-R_1(Q))=1.$$

Let $\dis g(\tH)=5+\frac{n-1}{2}$.
Then $\tH$ is either
$$2H+n\NI+(n+2)\NI+(n+4)\NI \mbox{ or } 2H+n\NI+(n+2)\NI+(n+6)\NI$$
$$ \mbox{ or } 2H+n\NI+(n+4)\NI+(n+6)\NI.$$
We will show that the first one is \dc.
Let $Q_1$ be a point of $C$ distinct from $P$.
There is a unique set consisting of three points $R_1$, $R_2$ and $R_3$ such that $4P-Q_1\sim R_1+R_2+R_3$.
Since we have $Q_1\not=P$ and $h^0(3P)=1$, we get $P\not\in\{R_1,R_2,R_3\}$.
We set $Q_2=R_1$ and $Q_3=R_2$.
Then we obtain
$$h^0(2P+Q_1+Q_2+Q_3)=2+h^0(4P-Q_1-Q_2-Q_3)=2+h^0(R_3)=3$$
and
$$h^0(P+Q_1+Q_2+Q_3)=1+h^0(5P-Q_1-Q_2-Q_3)=1+h^0(P+R_3)=2.$$
Moreover, we get
$h^0(Q_1+Q_2+Q_3)=1$.
In fact, we assume that $h^0(Q_1+Q_2+Q_3)=2$.
We have
$$2=h^0(2P+4P-Q_1-Q_2-Q_3)=h^0(2P+R_3).$$
If there exists a unique trigonal system $g_3^1$, then we get
$$6P\sim K\sim 2g_3^1\sim 4P+2R_3.$$
This contradicts $R_3\not=P$.
If there exist two trigonal systems $g_3^1$ and $h_3^1$,
Then we obtain
$$2=h^0(h_3^1)=h^0(K-g_3^1)=h^0(6P-2P-R_3)=h^0(4P-R_3).$$
This is a contradiction.
Next, we will prove that the second semigroup is \dc.
Let $Q_1$ be a point different from $P$.
For a general point $Q_2$ we have $h^0(5P-Q_1-Q_2)=1$.
Then for any point $Q_3$ we obtain
$$h^0(Q_1+Q_2+Q_3)=h^0(K-Q_1-Q_3-Q_2)=1.$$
On the other hand, we have
$5P-Q_1-Q_2\sim R_1+R_2+R_3$
where the set $\{R_1,R_2,R_3\}$ is uniquely determined.
Then we get $R_2\not=P$ or $R_3\not=P$, because of $h^0(3P)=1$ and $Q_1\not=P$.
Hence we may assume that $R_3\not=P$.
We set $Q_3=R_3$.
Then we obtain
$$h^0(P+Q_1+Q_2+Q_3)=1+h^0(5P-Q_1-Q_2-Q_3)=1+h^0(R_1+R_2)=2.$$
Moreover, we get
$$h^0(2P+Q_1+Q_2+Q_3)=2+h^0(4P-Q_1-Q_2-Q_3)=2,$$
because $Q_2$ is general.
On the other hand, we see that
$$h^0(3P+Q_1+Q_2+Q_3)=3+h^0(3P-Q_2-Q_1-Q_3)=3.$$
Finally, we will show that the last one is \dc.
Let $Q_1$, $Q_2$ and $Q_3$ be three points with $Q_1\not=P$.
Then we get $3P+Q_1+Q_2+Q_3\not\sim K$, because of $K\sim 6P$ and $Q_1\not=P$.
Hence, we get $h^0(3P+Q_1+Q_2+Q_3)=3$.
Moreover, we get
$$h^0(2P+Q_1+Q_2+Q_3)=2+h^0(4P-Q_1-Q_2-Q_3)=2$$
if $Q_2$ is general.
We also have
$$h^0(P+Q_1+Q_2+Q_3)=1+h^0(5P-Q_1-Q_2-Q_3)=1$$
if $Q_2$ and $Q_3$ are general.

Let $\dis g(\tH)=4+\frac{n-1}{2}$.
By Corollary \ref{lowest} $\tH$ is \dc.
\qed
\section{Numerical semigroups $\tH$ with $d_2(\tH)=\la 4,5,7\ra$
} 
\begin{theorem} 
We set $H=\la 4,5,7\ra$.
Let $\tH$ be a numerical semigroup with $d_2(\tH)=H$.
Assume that $g(\tH)\geqq 12$.
Then $\tH$ is \dc.
\end{theorem}
\pr
We take a pointed curve $(C,P)$ with $H(P)=\la 4,5,7\ra$.
In this case we have $K\sim 5P+R$ with $R\not=P$.
Let $n=n(\tH)\leqq 9$.
If $n=9$, then $\tH\supseteqq \la 8,9,10,14\ra$ whose genus is $12$.
Hence, we have $\tH=\la 8,9,10,14\ra$.
If we set $D=5P$, the same proof as in Theorem 2.2 of ~\cite{kom-ohb1} works well.
Thus, the semigroup $\tH$ is \dc.
If $n\leqq 7$, then we have $g(\tH)<12$.
Hence, we may assume that $n\geqq 11$.

Let $\dis g(\tH)=8+\frac{n-1}{2}$.
By Lemma \ref{ineq} and Proposition \ref{(n-1)/2} $\tH$ is of double covering type.

Let $\dis g(\tH)=7+\frac{n-1}{2}$.
We have the standard basis
$$S(2H+n\NI)=\{8,10,14,20,n,n+10,n+14,n+20\}.$$
Hence $\tH$ is either
$$2H+n\NI+(n+6)\NI \mbox{ or } 2H+n\NI+(n+12)\NI.$$
By Proposition \ref{smax-m} the second one is \dc.
Let $Q=R$.
Then we have
$$h^0(3P+Q)=1+h^0(2P+R-Q)=1+h^0(2P)=2 $$
and
$$h^0(2P+Q)=h^0(3P+R-Q)=h^0(3P)=1,$$
which implies that the first one is \dc.

Let $\dis g(\tH)=6+\frac{n-1}{2}$.
Then $\tH$ is either
$$2H+n\NI+(n+2)\NI \mbox{ or } 2H+n\NI+(n+4)\NI$$
$$\mbox{or }2H+n\NI+(n+6)\NI+(n+12)\NI.$$
We will show that the first semigroup is \dc.
Let $g_3^1$ be a trigonal system on $C$.
There is a unique set $\{Q_1,Q_2\}$ of points of $C$ with $g_3^1-P\sim Q_1+Q_2$.
Then we have
$$h^0(P+Q_1+Q_2)=2=h^0(Q_1+Q_2)+1.$$
Hence, it suffices to show that we may take $Q_1$ and $Q_2$ different from $P$.
Now we have $Q_1+Q_2\not=2P$.
Hence, we may assume that $Q_2\not=P$.
If $C$ has a unique trigonal system $g_3^1$.
Then we get
$$5P+R\sim K\sim 2g_3^1\sim 2P+2Q_1+2Q_2,$$
which implies that $3P+R\sim 2Q_1+2Q_2$.
If $Q_1=P$, then $P+R\sim 2Q_2$, which contradicts $Q_2\not= P$.
Assume that the curve $C$ has two trigonal systems $g_3^1$ and $h_3^1$ and $Q_1=P$.
Then we have
$$h_3^1\sim K-g_3^1\not\sim g_3^1\sim 2P+Q_2.$$
Let $h_3^1-P\sim S_1+S_2$ where the set $\{S_1,S_2\}$ is uniquely determined.
Then we may assume that $S_2\not=P$.
If $S_1=P$, then we get
$$5P+R\sim K\sim g_3^1+h_3^1\sim 4P+Q_2+S_2,$$
which implies that $P+R\sim Q_2+S_2$.
This contradicts $Q_2\not=P$ and $S_2\not=P$.
Thus, we replace $Q_1$ and $Q_2$ by $S_1$ and $S_2$ respectively.
For the third semigroup let $Q_1=R$ and let $Q_2$ be a general point.
Then we obtain
$$h^0(2P+Q_1+Q_2)=1+h^0(5P+R-2P-Q_1-Q_2)=1+h^0(3P-Q_2)=1$$
and
$$h^0(3P+Q_1+Q_2)=2+h^0(5P+R-3P-Q_1-Q_2)=2+h^0(2P-Q_2)=2.$$
Moreover, we have
$$h^0(5P+Q_1+Q_2)=4\mbox{ and }h^0(6P+Q_1+Q_2)=5.$$
Hence, the third one is \dc.
Let $C_0$ be a singular plane curve defined by the equation
\begin{equation}\label{f457}
a(yz-x^2)(y-2x+z)(y-2x-3z)z+bxy^2(y-z)^2=0.
\end{equation}
Using Bertini's theorem, for general $a,b\in k$, $C_0$ is a singular  plane curve with two singularities $A=(1,1,1)$ and $B=(-1,1,1)$, which are nodes.
In this case, let $C$ be the normalization of $C_0$.
Then the curve $C$ has genus $4$.
Moreover, let $P=(0,0,1) \in C$ be the origin.
Then we have
$\dim |3P|=-1+i(3P)$.
Let $C_2$ be a conic defined by $yz-x^2=0$.
Take a conic $C'$ in $i(3P)$, i.e., a conic $C'$ through $3P$, $A$ and $B$.
Then the intersection $C'.C_2$ contains $A$, $B$ and $3{P}$.
Hence we get $C'=C_2$, which implies that $\dim |3P|=0$.
Similarly, we obtain $\dim |4P|=1$ and $\dim |5P|=2$.
A conic in $i(6P)$ should be $C_2$.
But the intersection $C_0.C_2$ contains $5P$ but not $6P$.
Hence we get $i(6P)=0$, which implies that $\dim |6P|=2$.
Thus, we have $H(P)=\la 4,5,7\ra$.
Consider a line $L_0$ defined by the equation $y=0$.
Then the intersection $C_0.L_0$ is the divisor $2P+Q_1+Q_2+S$
where $P=(0,0,1)$, $\dis Q_1=(\frac{1}{2},0,1)$, $\dis Q_2=(-\frac{3}{2},0,1)$ and $S=(1,0,0)$.
Hence, we get 
$$\dim |2P+Q_1+Q_2|\geqq\dim |2P+Q_1+Q_2+S|-1$$
$$=1+h^0(K-2P-Q_1-Q_2-S)-1=1,$$
because a conic containing  $2P,Q_1,Q_2,S,A$ and $B$ should be $L_0L_1$
where $L_1$ is the line defined by the equation $y=z$.
Moreover, we have
$$\dim |P+Q_1+Q_2|=-1+h^0(K-P-Q_1-Q_2)=-1+1=0,$$
because a conic containing five points $P,Q_1,Q_2,A$ and $B$ should be $L_0L_1$.
Hence the second semigroup is \dc.

Let $\dis g(\tH)=5+\frac{n-1}{2}$.
Then $\tH$ is either
$$2H+n\NI+(n+4)\NI+(n+6)\NI \mbox{ or } 2H+n\NI+(n+2)\NI+(n+6)\NI$$
$$\mbox{or }2H+n\NI+(n+2)\NI+(n+4)\NI.$$
We will prove that the first semigroup is \dc.
Let $Q_1$ be a point of $C$ different from $P$ and $Q_i$ a general point for $i=2,3$.
Then we have
$$h^0(P+Q_1+Q_2+Q_3)=1+h^0(K-P-Q_1-Q_2-Q_3)=1,$$
$$h^0(2P+Q_1+Q_2+Q_3)=2\mbox{ and }h^0(3P+Q_1+Q_2+Q_3)=3.$$
Next, we will show that the second semigroup is \dc.
Let $Q_i\not=P$ for $i=1,2,3$ and $Q_j$ a general point for $j=1,2$.
Then $3P+Q_1+Q_2+Q_3\not\sim 5P+R\sim K$, because of $h^0(K-3P)=h^0(3P)=1$.
Hence, we get $h^0(3P+Q_1+Q_2+Q_3)=3$.
Moreover, we have
$$h^0(2P+Q_1+Q_2+Q_3)=2+h^0(K-2P-Q_1-Q_2-Q_3)=2+2-2=2,$$
because $Q_1$ and $Q_2$ are general.
Since $h^0(K-P-Q_1-Q_2)=1$, a divisor $S_1+S_2+S_3$ linearly equivalent to $K-P-Q_1-Q_2$ is uniquely determined. 
Then we have $S_1+S_2+S_3\not=3P$, because $Q_1$ and $Q_2$ are distinct from $P$.
Hence, we may assume that $S_3\not=P$.
We set $Q_3=S_3$.
Then we obtain
$$h^0(P+Q_1+Q_2+Q_3)=1+h^0(K-P-Q_1-Q_2-Q_3)=1+h^0(S_1+S_2)=2.$$
Since $Q_1$ and $Q_2$ are general, we may assume that $h^0(Q_1+Q_2+Q_3)=1$.
Lastly, we will see that the third one is \dc.
We use the curve $C$ which is the the normalization of the singular curve $C_0$ defined by (\ref{f457}).
Let $L_1$ and $L_2$ be the lines defined by the equations
$y=2x-z$ and $y=2x+3z$ respectively. 
The intersections $L_1.C_0$ and $L_2.C_0$ define distinct trigonal systems $g_3^1$ and $h_3^1$ on $C$.
We set $g_3^1-P\sim Q_1+Q_2$ and $h_3^1-P\sim Q_3+Q_4$.
In this case , we have $Q_1+Q_2\not\sim Q_3+Q_4$.
We may assume that $Q_1\not=P$ and $Q_3\not=P$, because $3\not\in H(P)$ implies that $g_3^1\not\sim 3P$ and $h_3^1\not\sim 3P$.
Moreover, we obtain
$$5P+R\sim K\sim g_3^1+h_3^1\sim 2P+Q_1+Q_2+Q_3+Q_4.$$
Since $Q_1\not=P$ and $Q_3\not=P$, the equalities $Q_2=Q_4=P$ do not hold.
Hence, we may assume that $Q_2\not=P$.
Then we get
$$h^0(2P+Q_1+Q_2+Q_3)=h^0(K-Q_4)=3$$
and
$$h^0(P+Q_1+Q_2+Q_3)=h^0(K-Q_4-P)=2.$$
Since $Q_i\not=P$ for $i=1,2,3$, we have $Q_1+Q_2+Q_3\not\sim g_3^1$ and $Q_1+Q_2+Q_3\not\sim h_3^1$.
Hence, we get $h^0(Q_1+Q_2+Q_3)=1$.

Let $\dis g(\tH)=4+\frac{n-1}{2}$.
By Corollary \ref{lowest} $\tH$ is \dc.
\qed
\section{Numerical semigroups $\tH$ with $d_2(\tH)=\la 4,6,7,9\ra$
} 
In this section we set $H=\la 4,6,7,9\ra$.
\begin{proposition} 
Let $\tH$ be a numerical semigroup with $d_2(\tH)=H$.
Assume that $g(\tH)\geqq 12$.
If $\dis g(\tH)=7+\fn$ where $n=n(\tH)$, then $\tH$ is of double covering type.
\end{proposition}
\pr
In this case, $\tH$ is either $2H+(n+10)\NI$ or $2H+(n+6)\NI$ or $2H+(n+4)\NI$.
Let $(C,P)$ be a pointed curve with $H(P)=\la 4,6,7,9\ra$.
Then we have $K\sim 4P+R_1+R_2$ with $R_i\not=P$ for $i=1,2$.
Let $Q$ be a point distinct from $P$, $R_1$ and $R_2$.
Then we get
$h^0(5P+Q)=3+h^0(R_1+R_2-P-Q)=3$ and
$h^0(4P+Q)=2+h^0(R_1+R_2-Q)=2$.
Hence, the first semigroup is \dc.

For general elements $a$ and $b\in k$, let $(C,P)$ be the pointed curve defined by the equation
$ax(y-x)^2(y+x)^2+b(yz^3-x^4)(y-z)=0$ with $P=(0,0,1)$.
Then the curve $C$ has only two nodes $A=(-1,1,1)$ and $B=(1,1,1)$ as singularities.
Let $T$ be the line defined by the equation $y=0$, which is the tangent line at $P$.
In fact, we have $T.C=4P+Q$ where $\dis Q=(-\frac{b}{a},0,1)$. 
We consider the line $L$ through $A$ and $B$, i.e., $y-z=0$.
Then we get $L.C=2A+2B+R$ where we set $R=(0,1,1)$.
Hence, the conic $y(y-z)$ cuts a canonical divisor $4P+Q+R$ on the normalization of $C$.
Since a conic containing $A,B,3P$ is uniquely determined, we get $h^0(3P)=h^0(K-3P)=1$, $h^0(4P)=1+h^0(K-4P)=2$ and $h^0(5P)=2+h^0(K-5P)=2$, which implies that $H(P)=\la 4,6,7,9\ra$.
On the other hand, we get $h^0(3P+Q)=h^0(K-P-R)=2$.
Moreover, we have
$h^0(2P+Q)=h^0(K-2P-Q)=1$, 
because a conic containing $A,B,2P,Q$ is uniquely determined, it is the conic $y(y-z)$.
Namely, the second semigroup is \dc.

For the third semigroup we use Proposition 17 in ~\cite{cop} as follows: There exists a curve $C$ of genus $4$ which is a triple covering of $\mathbb{P}^1$ with an ordinary ramification point $P$ whose Weierstrass semigroup is $\la 4,6,7,9\ra$.
Hence, we get $h^0(2P+Q)=2$ and $h^0(P+Q)=1$ for some point $Q$ distinct from $P$.
Namely, the third one is \dc.
\qed
\begin{proposition} 
Let $\tH$ be a numerical semigroup with $d_2(\tH)=H$.
Assume that $g(\tH)\geqq 12$.
If $\dis g(\tH)=6+\fn$ where $n=n(\tH)$, then $\tH$ is of double covering type.
\end{proposition}
\pr
Let $C_0$ be the singular plane curve defined by the equation
$$a(yz-x^2)(y-2x+z)(y-2x-3z)z+by^2(y-z)^2(x+2z)=0.$$
Using Bertini's theorem, for general $a$ and $b$ the curve $C_0$ has only two singularities $A=(1,1,1)$ and $B=(-1,1,1)$.
Let $C$ be the normalization of $C_0$.
Then the curve $C$ has genus $4$.
Let $P=(0,0,1)\in C$ be the origin.
Then we have
$\dim |3P|=-1+i(3P)$.
Let $C_2$ be the conic defined by $yz-x^2=0$.
Take a conic $C'$ in $i(3P)$.
Then the intersection $C'.C_2$ contains $A$, $B$ and $3P$.
Hence we get $C'=C_2$, which implies that $\dim |3P|=0$.
Similarly, we obtain $\dim |4P|=1$.
A conic in $i(5P)$ should be $C_2$.
But the intersection $C_0.C_2$ contains $4P$ but not $5P$.
Hence we get $i(5P)=0$, which implies that $\dim |5P|=1$.
Thus, we have $H(P)=\la 4,6,7,9\ra$.
Since we have
$$S(2H+n\NI)=\{8,12,14,18,n,n+12,n+14,n+18\},$$
The semigroup $\tH$ is either
$$2H+n\NI+(n+2)\NI\mbox{ or }2H+n\NI+(n+4)\NI+(n+6)\NI\mbox{ or}$$
$$2H+n\NI+(n+4)\NI+(n+10)\NI\mbox{ or }2H+n\NI+(n+6)\NI+(n+10)\NI.$$

We will prove that the first semigroup is \dc.
Consider the line $L$ defined by the equation $x-y=0$.
Then the intersection $C_0.L$ is the divisor $2A+P+Q_1+Q_2$
where $Q_1=(\alpha,\alpha,1)$ and $Q_2=(\beta,\beta,1)$ where $\alpha$ and $\beta$ are the roots of a quadratic equation
$bx^2+(2b-a)x-3a=0.$
Hence we get $|P+Q_1+Q_2|=g_3^1$, which implies that $2=h^0(P+Q_1+Q_2)=h^0(Q_1+Q_2)+1$.
Thus, the first one is \dc.

Next, we will show that the second semigroup is \dc.
The conic $C_2$ cuts out $C$ at $4P+R_1+R_2$ where $R_1$ and $R_2$ are distinct from $P$.
Hence, we get $\dim |4P+R_1+R_2|=\dim |K_C|=3$, $\dim |3P+R_1+R_2|=2$ and $\dim |2P+R_1+R_2|=1$.
Moreover, we have
$\dim |P+R_1+R_2|=0$,
because a conic through $A$, $B$, $P$, $R_1$ and $R_2$ should be $C_2$.
Thus, the second one is \dc.

We will prove that the third semigroup is \dc.
The line $y=0$ cuts out $C$ at $2P+S_1+S_2+R$ where $\dis S_1=(\frac{1}{2},0,1)$, $\dis S_2=(-\frac{3}{2},0,1)$ and $R=(1,0,0)$.
We have a unique conic containing $A$, $B$, $P$, $S_1$ and $S_2$ which has two irreducible components $y=0$ and the line through $A$ and $B$.
Hence, we get $\dim |P+S_1+S_2|=0$.
Moreover, we have $\dim |2P+S_1+S_2|=1$.
Since there is no conic cutting $3P+S_1+S_2$.
Hence, we get $\dim |iP+S_1+S_2|=i-2$ for $i=3,4,5$.
Thus, the third one is \dc.

Lastly, we will show that the fourth semigroup is \dc.
Take a point $Q$ in the intersection of $yz-x^2=0$ and $C$, which is distinct from $P$, $A$ and $B$.
Let $S$ be a point not in in the intersection of $yz-x^2=0$ and $C$.
Then there is no conic containing $2P+Q+S+A+B$.
Hence, we get $\dim |jP+Q+S|=j-2$ for all $j=2,3,4,5$.
Thus, the fourth one is \dc.
\qed
\begin{proposition} 
Let $\tH$ be a numerical semigroup with $d_2(\tH)=H$.
Assume that $g(\tH)\geqq 12$.
If $\dis g(\tH)=5+\frac{n-1}{2}$ where $n=n(\tH)$, then $\tH$ is of double covering type.
\end{proposition}
\pr
Let $(C,P)$ be a pointed curve with $H(P)=\la 4,6,7,9\ra$.
Then we have $K\sim 4P+S_1+S_2$ where $S_1$ and $S_2$ are different from $P$.
In this case the semigroup $\tH$ is either $2H+\la n+4,n+6,n+10\ra$ or $2H+\la n+2,n+6\ra$ or $2H+\la n+2,n+4\ra$.

We will prove that the first semigroup is \dc.
Let $Q_1\not=P$ and let $Q_2$ and $Q_3$ be general points of $C$.
Then we have
$h^0(K-P-Q_1-Q_2-Q_3)=0$.
Hence, we get
$$5=h^0(5P+Q_1+Q_2+Q_3)=h^0(4P+Q_1+Q_2+Q_3)+1 \mbox{ and }$$
$$h^0(3P+Q_1+Q_2+Q_3)=3=h^0(2P+Q_1+Q_2+Q_3)+1=h^0(P+Q_1+Q_2+Q_3)+2.$$
Hence, the first one is \dc.

We will show that the second semigroup is \dc.
Let $Q_1$ and $Q_2$ be general points of $C$.
Then we get
$$h^0(3P+Q_1+Q_2+Q_3)=3\mbox{ and }h^0(2P+Q_1+Q_2+Q_3)=2.$$
Let us consider
$$h^0(P+Q_1+Q_2+Q_3)=1+h^0(K-P-Q_1-Q_2-Q_3).$$
We have
$K-P-Q_1-Q_2\sim R_1+R_2+R_3$ where the set $\{R_1,R_2,R_3\}$ is uniquely determined.
Moreover, $R_1+R_2+R_3\not=3P$, because $R_1+R_2+R_3=3P$ implies that $Q_1+Q_2=S_1+S_2$, which is a contradiction.
We may assume that $R_3\not=P$.
We set $Q_3=R_3$.
Then we get
$$h^0(P+Q_1+Q_2+Q_3)=1+h^0(R_1+R_2+R_3-Q_3)=1+h^0(R_1+R_2)=2.$$
Moreover, $h^0(Q_1+Q_2+Q_3)=1$, because $Q_1$ and $Q_2$ are general points.
Hence, the second one is \dc.

Finally, we will prove that the third semigroup is \dc.
We may assume that $C$ is a general curve with an ordinary Weierstrass point $P$, i.e., $H(P)=\{4,6,7,9\}$.
Then $C$ has two trigonal systems $g_3^1$ and $h_3^1$ with $g_3^1\not\sim h_3^1$.
Let $Q_1$ be a general point of $C$.
Then $K-2P-Q_1\sim Q_2+Q_3+Q_4$ where the set $\{Q_2,Q_3,Q_4\}$ is uniquely determined.
We note that
$$h^0(2P+Q_i+Q_j+Q_k)=2+h^0(K-2P-Q_i-Q_j-Q_k)=2+h^0(Q_l)=3$$
for $\{i,j,k,l\}=\{1,2,3,4\}$.
Moreover, we have
$$h^0(P+Q_i+Q_j+Q_k)=1+h^0(K-P-Q_i-Q_j-Q_k)=1+h^0(P+Q_l)=2$$
for $\{i,j,k,l\}=\{1,2,3,4\}$.
Hence, it suffices to show that there are distinct $i,j,k$ such that $h^0(Q_i+Q_j+Q_k)=1$ where three points $Q_i$, $Q_j$ and $Q_k$ are different from $P$.
We may assume that $Q_2$ and $Q_3$ are distinct from $P$, because if not, then 
$4P+S_1+S_2\sim 2P+Q_1+Q_2+Q_3+Q_4$ implies that $S_1+S_2\sim Q_1+Q_i$ for some $i$, hence a contradiction.
We assume that $h^0(Q_1+Q_2+Q_3)=2$.
Then we have
$$g_3^1\sim Q_1+Q_2+Q_3 \mbox{ and } h_3^1\sim K-(Q_1+Q_2+Q_3)\sim 2P+Q_4.$$
Hence, we get $Q_4\not=P$.
Here we assume that
$$h^0(Q_1+Q_2+Q_3)=h^0(Q_1+Q_2+Q_4)=h^0(Q_1+Q_3+Q_4)=h^0(Q_2+Q_3+Q_4)=2.$$
Using $h^0(Q_1+Q_2+Q_4)=2$, we get $Q_1+Q_2+Q_4\sim Q_1+Q_2+Q_3$ or $Q_1+Q_2+Q_4\sim 2P+Q_4$.
Since $Q_1$ is general, we have $Q_4\sim Q_3$, hence $Q_4=Q_3$.
By $h^0(Q_1+Q_3+Q_4)=2$ we get $Q_1+Q_3+Q_4\sim Q_1+Q_2+Q_3$ or $Q_1+Q_3+Q_4\sim 2P+Q_4$.
Hence, we have $Q_4=Q_2$.
Thus we have
$$2=h^0(Q_2+Q_3+Q_4)=h^0(3Q_4),$$
which implies that 
$3Q_4\sim Q_1+2Q_4$ or $3Q_4\sim 2P+Q_4$.
In view of $Q_4\not=P$ we have $Q_4=Q_1$, which implies that $h^0(3Q_1)=2$.
Since $Q_1$ is general, this is a contradiction.
Hence, the third one is \dc.
\qed
\begin{theorem} 
Let $\tH$ be a numerical semigroup with $d_2(\tH)=H$.
Assume that $g(\tH)\geqq 12$.
Then $\tH$ is \dc.
\end{theorem}
\pr
We set $n=n(\tH)$.
If $n\leqq 7$, then $\tH\supseteqq \la 7,8,12,18\ra$ whose genus is $11$.
Thus, we may assume that $n\geqq 9$.
Hence, by Lemma \ref{ineq} we may assume that $\dis 4+\frac{n-1}{2}\leqq g(\tH)\leqq 8+\frac{n-1}{2}$.
By Propositions \ref{(n-1)/2}, 7.1, 7.2, 7.3 and Corollary \ref{lowest} the numerical semigroup $\tH$ is \dc.
\qed
\section{Numerical semigroups $\tH$ with $d_2(\tH)=\la 5,6,7,8,9\ra$
} 
In this section we set $H=\la 5,6,7,8,9\ra$.
\begin{proposition} 
Let $\tH$ be a numerical semigroup with $d_2(\tH)=H$.
Assume that $g(\tH)\geqq 12$ and $n=n(\tH)\geqq 11$.
If $\dis g(\tH)=7+\frac{n-1}{2}$, then $\tH$ is \dc.
\end{proposition}
\pr
The semigroup $\tH$ is either 
$$2H+(n+8)\NI\mbox{ or }2H+(n+6)\NI\mbox{ or }
$$
$$2H+(n+4)\NI\mbox{ or }2H+(n+2)\NI.
$$

For the first semigroup by Proposition \ref{smax-m} it is \dc.

For the last semigroup let $(C,P)$ be a pointed hyperelliptic curve with an ordinary point $P$.
Then there is a unique point $Q\in C$, distinct from $P$, such that $h^0(P+Q)=2=h^0(Q)+1$.
Hence, the last one is \dc.

We will prove that the second semigroup is \dc.
Let $C$ be a curve of genus $4$ with a unique trigonal system $g_3^1$ and $P$ its ordinary point.
Then we have $K\sim 3P+R_1+R_2+R_3$ with $R_i\not=P$ for each $i$.
We set $Q=R_3$.
Then we obtain
$h^0(3P+Q)=1+h^0(R_1+R_2)=2$.
Moreover, we have $h^0(2P+Q)=h^0(P+R_1+R_2)$.
Assume that $h^0(P+R_1+R_2)=2$.
Then we get $g_3^1\sim P+R_1+R_2$, which implies that
$$3P+R_1+R_2+R_3\sim K\sim 2g_3^1\sim 2P+2R_1+2R_2.$$
Hence, we obtain $P+R_3\sim R_1+R_2$, which contradicts $R_1\not= P$ and $R_2\not=P$.
Thus, we get $h^0(2P+Q)=h^0(P+R_1+R_2)=1$.
Therefore, the second one is \dc.

Lastly we will show that the third semigroup is \dc.
Let $C_0$ be the singular plane curve defined by the equation (\ref{f457}).
By Bertini's theorem, for general $a$ and $b$ the curve $C_0$ has only two singularities $A=(1,1,1)$ and $B=(-1,1,1)$.
Let $C$ be the normalization of $C_0$.
Then the curve $C$ has genus $4$.
We set $\dis P=(\frac{1}{2},0,1)$ and $S=(0,-1,1)$.
Then the intersection $C_0.(y-2x+z)$ is $2P+S+2A$, which determines a trigonal system $g_3^1\sim 2P+S$.
Moreover, we have $\dim |4P|=i(4P)=0$, i.e., $P$ is an ordinary point.
In fact, let $C_2$ be a conic such that the intersection $C_2.C$ contains $4P+A+B$.
Then the intersection $C_2.(y-2x+z)$ contains $2P+A$.
Hence, $C_2$ is reducible, so $C_2=(y-2x+z)\cdot l$ for a line $l$ such that the intersection $l.C_0$ contains $2P+B$.
This is a contradiction, because there are two different tangent lines at $P$.
Hence, we get $i(4P)=0$.
We set $Q=S$.
Then we get
$$h^0(P+Q)=1\mbox{ and }h^0(2P+Q)=h^0(g_3^1)=2.$$
Thus, the third one is \dc.
\qed
\begin{proposition} \label{56789+6} 
Let $\tH$ be a numerical semigroup with $d_2(\tH)=H$.
Assume that $g(\tH)\geqq 12$ and $n=n(\tH)\geqq 11$.
If $\dis g(\tH)=6+\frac{n-1}{2}$, then $\tH$ is \dc.
\end{proposition}
\pr
The semigroup $\tH$ is either 
$$2H+(n+6)\NI+(n+8)\NI\mbox{ or }2H+(n+4)\NI+(n+8)\NI
$$
$$\mbox{or }2H+(n+4)\NI+(n+6)\NI\mbox{ or }2H+(n+2)\NI+(n+8)\NI$$
$$\mbox{or }2H+(n+2)\NI+(n+6)\NI\mbox{ or }2H+(n+2)\NI+(n+4)\NI.
$$

We will prove that the first semigroup is \dc.
We take a pointed curve $(C,P)$ with an ordinary point $P$.
Then we have $K\sim 3P+R_1+R_2+R_3$ where the divisor $R_1+R_2+R_3$ is uniquely determined.
Let $Q_1$ and $Q_2$ be general points.
Then we have
$$h^0(4P+Q_1+Q_2)=3\mbox{ and }h^0(3P+Q_1+Q_2)=2+h^0(K-2P-Q_1-Q_2-P)=2.$$
Moreover, we get
$$h^0(2P+Q_1+Q_2)=1+h^0(K-2P-Q_1-Q_2)=1.$$

We will see that the second one is \dc.
Let $C$ be a curve of genus $4$ with a unique trigonal system $g_3^1$.
We take an ordinary point $P$ of $C$.
We note that $K\sim 3P+R_1+R_2+R_3$ with $R_i\not=P$ for all $i$.
Let $Q_1$ be a general point and $K-2P-Q_1\sim Q_2+Q_3+Q_4$.
Then we get
$$h^0(4P+Q_1+Q_i)=3+h^0(K-4P-Q_1-Q_i)=3$$
and
$$h^0(3P+Q_1+Q_i)=2+h^0(K-3P-Q_1-Q_i)=2+h^0(R_1+R_2+R_3-Q_1-Q_i)=2.$$
Now we have
$$h^0(2P+Q_1+Q_i)=1+h^0(K-2P-Q_1-Q_i)=1+h^0(Q_j+Q_k)=2$$
where $\{i,j,k\}=\{2,3,4\}$.
We assume that
$h^0(P+Q_1+Q_2)=2$.
Then we get
$$2P+Q_1+Q_2+Q_3+Q_4\sim K\sim 2g_3^1\sim 2P+2Q_1+2Q_2,$$
which implies that
$Q_1+Q_2\sim Q_3+Q_4$.
Hence, we may assume that $Q_1=Q_3$ and $Q_2=Q_4$.
Assume that $h^0(P+Q_1+Q_3)=2$.
Then we get
$$2P+2Q_1+2Q_2\sim K\sim 2P+2Q_1+2Q_3\sim 2P+4Q_1,$$
which implies that $2Q_2\sim 2Q_1$.
Hence, we have $Q_1=Q_2$, so $Q_1=Q_2=Q_3=Q_4$.
Thus, we get $K\sim 2P+4Q_1$.
Hence, we have
$h^0(4Q_1)=h^0(K-2P)=2,$
which contradicts the generality of $Q_1$.
Thus, we get
$h^0(P+Q_1+Q_2)=1$.

We will prove that the third semigroup is \dc. 
Let $C$, $P$ and $R_i$'s be as in the above case.
We set $Q_1=R_1$ and $Q_2=R_2$.
We have
$$h^0(3P+Q_1+Q_2)=h^0(K-R_3)=3$$
and
$$h^0(2P+Q_1+Q_2)=h^0(K-R_3-P)=2.$$
We assume that
$h^0(P+Q_1+Q_2)=2$.
Then we have
$$3P+R_1+R_2+R_3\sim K\sim 2P+2Q_1+2Q_2,$$
which implies that $P+R_3\sim R_1+R_2$.
This contradicts $P\not=R_i$ for $i=1,2$.
Thus, we obtain $h^0(P+Q_1+Q_2)=1$.

We will prove that the fourth semigroup is \dc.
Let $C$, $P$ and $R_i$'s be as in the above case.
Let $g_3^1-P\sim Q_1+Q_2$.
Then we have $Q_i\not=P$ for $i=1,2$.
In fact, assume that $Q_1=P$.
Then we obtain $$3P+R_1+R_2+R_3\sim K\sim 2g_3^1\sim 4P+2Q_2,$$
which implies that $R_1+R_2+R_3\sim P+2Q_2$.
This contradicts $R_i\not=P$, all $i$.
Hence, we get
$$h^0(Q_1+Q_2)=1\mbox{ and }h^0(P+Q_1+Q_2)=h^0(g_3^1)=2.$$
Moreover, we have
$$h^0(3P+Q_1+Q_2)=2+h^0(2g_3^1-3P-Q_1-Q_2)=2+h^0(Q_1+Q_2-P)=2$$
and
$$h^0(4P+Q_1+Q_2)=3+h^0(2g_3^1-4P-Q_1-Q_2)=3.$$

We will prove that the fifth semigroup is \dc.
We use the notations in the proof where the third semigroup in Proposition 8.1 is of double covering type.
Let $K\sim 3P+R_1+R_2+R_3$ with $R_i\not=P$, all $i$.
Then the curve $C$ defined by the equation (\ref{f457}) has two distinct trigonal systems $g_3^1$ and $h_3^1$.
If not, we have $4P+2S\sim 2g_3^1\sim K\sim 3P+R_1+R_2+R_3$, which implies that $P+2S\sim R_1+R_2+R_3$.
This contradicts $R_i\not=P$, all $i$.
Let $h_3^1-P\sim S'+S''$.
Then we get
$$3P+S+S'+S''\sim g_3^1+h_3^1\sim K\sim 3P+R_1+R_2+R_3,$$
which implies that $S+S'+S''=R_1+R_2+R_3$.
Hence, we may assume that $S=R_1$, $S'=R_2$ and $S''=R_3$, so we have
$g_3^1\sim 2P+R_1$ and $h_3^1\sim P+R_2+R_3$.
We set $Q_1=R_2$ and $Q_2=R_3$.
Then we have
$$h^0(Q_1+Q_2)=1\mbox{ and }h^0(P+Q_1+Q_2)=h^0(h_3^1)=2.$$
Moreover, we get
$$h^0(2P+Q_1+Q_2)=1+h^0(K-2P-Q_1-Q_2)=1+h^0(P+R_1)=2$$
and
$$h^0(3P+Q_1+Q_2)=2+h^0(K-3P-Q_1-Q_2)=2+h^0(R_1)=3.$$

We will show that the last semigroup is \dc.
Let $C$ be a hyperelliptic curve of genus $4$.
Then there is a unique point $Q$ distinct from $P$ with $P+Q\sim g_2^1$.
We set $Q_1=Q_2=Q$.
Then we have 
$$h^0(Q_1+Q_2)=h^0(2Q)=1\mbox{ and }h^0(P+Q_1+Q_2)=h^0(P+Q+Q)=2.$$
Moreover, we obtain
$$h^0(2P+Q_1+Q_2)=1+h^0(3g_2^1-2P-2Q)=1+h^0(g_2^1)=3.$$
\qed
\begin{proposition} 
Let $\tH$ be a numerical semigroup with $d_2(\tH)=H$.
Assume that $g(\tH)\geqq 12$ and $n=n(\tH)\geqq 11$.
If $\dis g(\tH)=5+\frac{n-1}{2}$, then $\tH$ is of double covering type.
\end{proposition}
\pr
The semigroup $\tH$ is either 
$$2H+(n+4)\NI+(n+6)\NI+(n+8)\NI\mbox{ or }2H+(n+2)\NI+(n+6)\NI+(n+8)\NI\mbox{ or }
$$
$$2H+(n+2)\NI+(n+4)\NI+(n+8)\NI\mbox{ or }2H+(n+2)\NI+(n+4)\NI+(n+6)\NI.$$
We take a pointed curve $(C,P)$ with an ordinary point $P$.
Then we have $K\sim 3P+R_1+R_2+R_3$ where the divisor $R_1+R_2+R_3$ is uniquely determined.

We will show that the first semigroup is \dc.
Let $Q_1$, $Q_2$ and $Q_3$ be general points of $C$.
Then we have
$h^0(4P+Q_1+Q_2+Q_3)=4$ and $h^0(3P+Q_1+Q_2+Q_3)=3.$
Moreover, we get
$$h^0(2P+Q_1+Q_2+Q_3)=2+h^0(K-2P-Q_1-Q_2-Q_3)=2$$
and
$$h^0(P+Q_1+Q_2+Q_3)=1+h^0(K-P-Q_1-Q_2-Q_3)=1.$$

Next, we will prove that the second semigroup is \dc.
Let $Q_1$ and $Q_2$ be general points and $Q_3$ a point distinct from $P$.
We have $h^0(4P+Q_1+Q_2+Q_3)=4$.
Since $Q_1$ and $Q_2$ are general, we get
$$h^0(3P+Q_1+Q_2+Q_3)=3 \mbox{ and }
h^0(2P+Q_1+Q_2+Q_3)=2.$$
Let $K-P-Q_1-Q_2\sim S_1+S_2+S_3$ for some points $S_1$, $S_2$ and $S_3$.
If $S_1+S_2+S_3=3P$, then we get $R_1+R_2+R_3\sim P+Q_1+Q_2$, which contradicts $R_i\not=P$ for all $i$.
Hence, we may assume that $S_3\not=P$.
We set $Q_3=S_3$.
Then we get
$$h^0(P+Q_1+Q_2+Q_3)=1+h^0(K-P-Q_1-Q_2-Q_3)=1+h^0(S_1+S_2)=2.$$
Since $Q_1$ and $Q_2$ are general, we have $h^0(Q_1+Q_2+Q_3)=1$.

We will prove that the third semigroup is \dc.
Assume that $C$ has a unique trigonal system $g_3^1$.
Let $Q_1$ be a general point, and $Q_2$ and $Q_3$ two points distinct from $P$.
Then we have
$$h^0(4P+Q_1+Q_2+Q_3)=4\mbox{ and }h^0(3P+Q_1+Q_2+Q_3)=3.$$
Let $K-2P-Q_1\sim S_1+S_2+S_3$.
Then we have $S_i\not=P$.
In fact, let $S_3=P$.
Then we get
$$R_1+R_2+R_3\sim Q_1+S_1+S_2,$$
which contradicts the generality of $Q_1$.
We set
$Q_2=S_2$ and $Q_3=S_3$.
Then we have
$$h^0(2P+Q_1+Q_2+Q_3)=2+h^0(K-2P-Q_1-Q_2-Q_3)=2+h^0(S_1)=3.$$
Moreover, we see that
$$h^0(P+Q_1+Q_2+Q_3)=1+h^0(K-P-Q_1-Q_2-Q_3)=1+h^0(P+S_1)=2$$
and
$$h^0(Q_1+Q_2+Q_3)=h^0(K-Q_1-Q_2-Q_3)=h^0(2P+S_1).$$
We assume that $h^0(2P+S_1)=2$.
Since $g_3^1$ is a unique trigonal system which is linearly equivalent to $2P+S_1$, we have
$$3P+R_1+R_2+R_3\sim K\sim 4P+2S_1,$$
which induces that $R_1+R_2+R_3\sim P+2S_1$.
This contradicts $R_i\not=P$ for all $i$.
Hence, we get $h^0(Q_1+Q_2+Q_3)=1$.

Finally, we will see that the last semigroup is \dc.
We set $Q_i=R_i$ for $i=1,2,3$.
Then we have
$$h^0(3P+Q_1+Q_2+Q_3)=h^0(K)=4\mbox{, }h^0(2P+Q_1+Q_2+Q_3)=h^0(K-P)=3,$$
$$h^0(P+Q_1+Q_2+Q_3)=h^0(K-2P)=2\mbox{ and }h^0(Q_1+Q_2+Q_3)=h^0(R_1+R_2+R_3)=1.$$
\qed
\begin{theorem} 
Let $\tH$ be a numerical semigroup with $d_2(\tH)=H$.
Assume that $g(\tH)\geqq 12$.
Then $\tH$ is \dc.
\end{theorem}
\pr
Let $n=n(\tH)\leqq 9$.
If $n=9$, then we get $\tH\supseteqq \la 9,10,12,14,16\ra$ whose genus is $12$.
Let $(C,P)$ be a pointed curve with an ordinary point $P$.
We set $D=5P$.
The same proof as in Theorem 2.2 of ~\cite{kom-ohb1} works well.
Thus, the semigroup $\tH$ is \dc.
If $n\leqq 7$, then we have $g(\tH)<12$.
Hence, we may assume that $n\geqq 11$.
If $\dis g(\tH)=8+\frac{n-1}{2}$, $\dis 7+\frac{n-1}{2}$, $\dis 6+\frac{n-1}{2}$, $\dis 5+\frac{n-1}{2}$, $\dis 4+\frac{n-1}{2}$, then $\tH$ is \dc \ by Propositions \ref{(n-1)/2}, 8.1, 8.2, 8.3, Corollary \ref{lowest} respectively.
\qed

\end{document}